\documentclass{amsart}
\usepackage[T1]{fontenc}
\usepackage{pinlabel}
\usepackage{amsmath, amsthm, amssymb, amscd, mathrsfs, eucal, epsfig, color}
\usepackage{hyperref}

\input xy
\xyoption{all}

\begin{document}

\newtheorem{theorem}{Theorem}[section]
\newtheorem{lemma}[theorem]{Lemma}
\newtheorem{corollary}[theorem]{Corollary}
\newtheorem{conjecture}[theorem]{Conjecture}
\newtheorem{proposition}[theorem]{Proposition}
\newtheorem{question}[theorem]{Question}
\newtheorem*{answer}{Answer}
\newtheorem{problem}[theorem]{Problem}
\newtheorem*{main_theorem}{Main Theorem}
\newtheorem*{claim}{Claim}
\newtheorem*{criterion}{Criterion}
\theoremstyle{definition}
\newtheorem{definition}[theorem]{Definition}
\newtheorem{construction}[theorem]{Construction}
\newtheorem{notation}[theorem]{Notation}
\newtheorem{convention}[theorem]{Convention}
\newtheorem*{warning}{Warning}
\newtheorem*{assumption}{Simplifying Assumptions}

\theoremstyle{remark}
\newtheorem{remark}[theorem]{Remark}
\newtheorem*{apology}{Apology}
\newtheorem{historical_remark}[theorem]{Historical Remark}
\newtheorem{example}[theorem]{Example}
\newtheorem{scholium}[theorem]{Scholium}
\newtheorem*{case}{Case}

\def\T{\mathcal T}
\def\M{\mathcal M}
\def\N{\mathbb N}
\def\R{\mathbb R}
\def\Z{\mathbb Z}
\def\D{\mathbb D}
\def\H{\mathbb H}
\def\C{\mathbb C}
\def\Chat{{\mathbb CP}^1}

\newcommand{\marginal}[1]{\marginpar{\tiny #1}}

\title{Surgery sequences and self-similarity of the Mandelbrot set}
\author{Danny Calegari}
\address{Department of Mathematics \\ University of Chicago \\
Chicago, Illinois, 60637}
\email{dannyc@math.uchicago.edu}

\date{\today}
\dedicatory{for Tan Lei}

\begin{abstract}
We introduce an analog in the context of rational maps of the idea of
{\em hyperbolic Dehn surgery} from the theory of Kleinian groups. A {\em surgery sequence}
is a sequence of postcritically finite maps limiting (in a precise manner) to
a postcritically finite map with at least one strictly preperiodic critical orbit. 
As an application of this idea we give a new and
elementary proof of Tan Lei's theorem on the asymptotic self-similarity of Julia sets and
the Mandelbrot set at Misiurewicz points.
\end{abstract}

\maketitle

\setcounter{tocdepth}{1}
\tableofcontents

\section{Introduction}

{\em Sullivan's Dictionary} is a framework that seeks to unify two subjects at the heart
of one-dimensional holomorphic dynamics: Kleinian groups; and the iteration of rational maps.
Sullivan introduced the idea of this dictionary in \cite{Sullivan}, and supplied a number
of key entries. It is the purpose of this paper to propose a new entry for this dictionary
--- between {\em hyperbolic Dehn surgery} for cusped hyperbolic 3-manifolds and {\em surgery
sequences} of rational maps --- and to use this analogy to give a new, elementary proof
of Tan Lei's famous theorem on the asymptotic self-similarity of Julia sets and
the Mandelbrot set at Misiurewicz points.

Sullivan's original dictionary contains twenty-four entries; the last entry is a correspondence
between cocompact Kleinian groups and postcritically finite rational maps. One might
reasonably broaden this correspondence from the class of cocompact Kleinian groups to the
finite covolume groups. The idea of this dictionary entry is that in either case,
topological data (an irreducible compact 3-manifold with torus boundary components; an
equivalence class of postcritically finite branched self-covering map of the 2-sphere)
may be `geometrized' by a rigid holomorphic dynamical system, unless a purely 
topological obstruction exists.

One reason to include finite covolume Kleinian groups in the picture is that these groups arise
as {\em limits} of the cocompact ones. A rank 2 parabolic subgroup $H$ of a Kleinian group $G$
corresponds to a {\em toral cusp} in the 3-manifold quotient $M$. Dehn filling a cusp of $M$
with a long slope gives rise to a new 3-manifold $M'$, which is the quotient of $\H^3$
by a new Kleinian group $G'$, which has a rank 1 loxodromic subgroup $H'$ `in place' of the
rank 2 parabolic subgroup $H$ of $G$. For suitable choice of filling slope,
the geometry and topology of $M'$ (resp. the dynamics of $G'$
on the Riemann sphere) will approximate arbitrarily closely the geometry and topology 
of $M$ away from the cusp (resp. the dynamics of $G$ away from the fixed points 
of the conjugates of $H$).

Now let's move to the other side of the dictionary. What is the analog of Dehn filling in the
world of rational maps?
For a postcritically finite rational map $f$ there is no really good analog of the quotient 
3-manifold $M$ on the Kleinian group side, 
and any sort of approximation must take place in the dynamics on the Riemann
sphere. We propose the following informal analogy; for a precise definition see
Definition~\ref{definition:surgery_sequence}. 
We start with the data of a postcritically finite map $f$ with
at least one critical point $c$ whose forward orbit contains a repelling cycle $O$ disjoint from $c$.
The cycle $O$ is the analog of the `cusp' which is to be deformed. We must also make a
choice of an infinite backward orbit $T$ of $c$ (called the `tail') accumulating only on $O$. A
{\em surgery sequence} is then a sequence of postcritically finite maps $f_n$ so that
$f_n \to f$ as maps, and so that the postcritical set of $f_n$ converges (in the Hausdorff
topology) to the {\em union} of the postcritical set of $f$ with $T$.

Different choices of tail $T$ accumulating on $O$ give rise to different surgery
sequences for a fixed $f$; we may describe the structure of the set of all such surgery sequences,
and this description lets us
recover the (asymptotic) geometry of the Julia set $J(f)$ of $f$ near the orbit $O$, and
in the special case that $f$ has degree 2 and corresponds to a Misiurewicz point in the 
Mandelbrot set, the (asymptotic) geometry of the Mandelbrot set near $f$. 

\section{Surgery sequences}

\subsection{Definition}\label{definition_subsection}

Let $f$ be a rational map. We denote the {\em critical set} of $f$ (i.e.\/ the set of critical
points of $f$) by $C(f)$. For each $c\in C(f)$ let $P(c):=\cup_{n>0} f^n(c)$ denote the
forward orbit of $c$, and let $P(f):=\cup_{c\in C(f)} P(c)$ denote the {\em postcritical set}.
The map $f$ is {\em postcritically finite} (pcf) if $P(f)$ is finite.

Let $f$ be a pcf map with $c\in C(f)$ and let $O \subset P(c)-C(f)$ be a  
periodic orbit with multiplier $\mu:=\prod_{x\in O} f'(x)$. Let's suppose $O$ is a repelling
orbit, i.e.\/ $|\mu|>1$.

Since $O$ is repelling, there is a neighborhood $U$ of $O$ and a unique branch $g$ of
$f^{-1}$ with $g:U \to U$, so that for every point $x\in U$ the sequence $g^n(x)$ accumulates
on $O$. A {\em tail} $T$ for $c,O$ is an infinite sequence $c_{-n}$ for $n\in \Z^+$ for which
\begin{enumerate}
\item{$f(c_{-1})=c$;}
\item{$f(c_{-n})=c_{1-n}$ for $n>1$; and}
\item{for $n$ sufficiently large, $c_{-n}\in U$ and $g(c_{-n})=c_{-1-n}$.}
\end{enumerate}

\begin{definition}\label{definition:surgery_sequence}
Let $f$ be pcf. Let $c$ be critical for $f$ and let $O \subset P(c)-C(f)$ be a repelling periodic
orbit for $f$. A {\em surgery sequence} for $c,O$ is a sequence of
pcf maps $f_n$ for which there is a tail $T$ for $c,O$ such that
\begin{enumerate}
\item{the $f_n$ converge to $f$ as rational maps;}
\item{the cardinality of $P(f_n)$ is $A+|O|n$ for some integer $A$; and}
\item{the sets $P(f_n)$ converge in the Hausdorff topology to $T\cup P(f)$.}
\end{enumerate}
\end{definition}

The idea of a surgery sequence is to approximate the `orbit' $T \cup c \cup P(c)$ of $f$ 
by some finite periodic critical orbit of $f_n$. 

We shall consider two surgery sequences $f_n$, $g_n$ to be {\em isomorphic} either
\begin{enumerate}
\item{if they are conjugate by a  (convergent) family of M\"obius transformations $\phi_n$
--- i.e.\/ $f_n = \phi_n g_n \phi_n^{-1}$ with $\phi_n \to \text{id}$; or}
\item{if there is an integer $m$ with $f_n = g_{n+m}$ for all $n$ (where defined).}
\end{enumerate}

In practice the first kind of ambiguity is eliminated by normalizing 
our surgery sequences somehow, either by fixing a small number of specific critical 
points and their images, or by fixing some of the coefficients of the rational maps $f_n$.

\subsection{An example}\label{2_example}

\begin{figure}[htpb]
\centering
\includegraphics[scale=0.4]{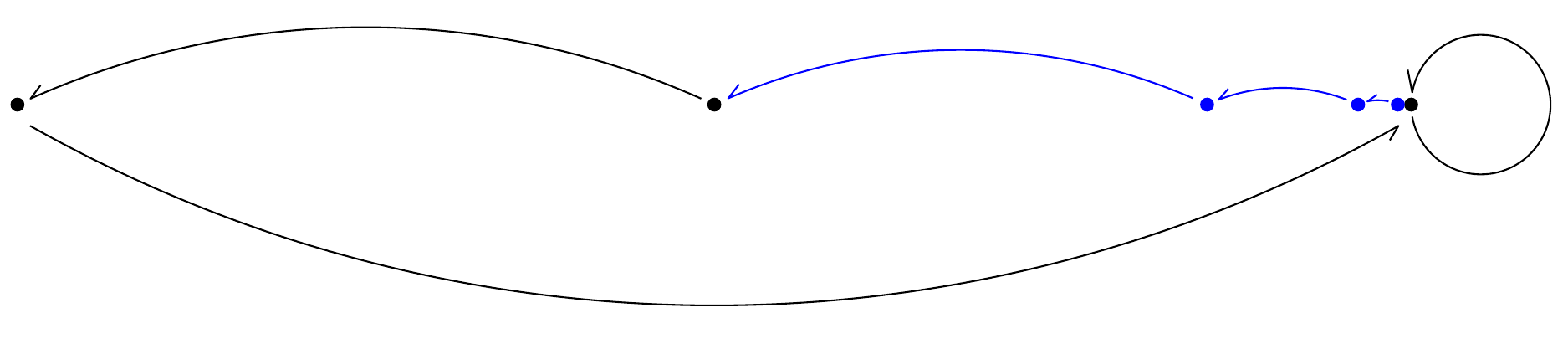}
\caption{The forward orbit $P(0)$ for the map $f:z \to z^2 - 2$ is in black, and the tail $T$ is in blue
(only three elements of $T$ are shown in the figure).}
\label{f_orbit}
\end{figure}

We consider the simplest nontrivial example, a surgery sequence for the pcf map 
$f:z \to z^2 - 2$. Since $f$ is a polynomial, the critical point $\infty$ is completely invariant. 
The only other critical point is $0$, whose orbit is $0 \to -2 \to 2$ thereafter
left fixed by $f$. Let $O$ consist solely of the point $p=2$; the multiplier $\mu = 4$ so
this is indeed a repelling fixed point; and as tail $T$ we choose the
sequence $c_{-1}=\sqrt{2}$ and $c_{-n} = \sqrt{c_{1-n}+2}$, where by convention, we take the 
`square root' symbol to denote the unique non-negative real square root of a non-negative real number.
Vi\`ete's formula \cite{Viete} implies that 
$$\lim_{n\to\infty} (c_{-n}-p)\cdot \mu^n = \lim_{n\to\infty} (c_{-n}-2)\cdot 4^n = -\pi^2/4 \sim -2.4675011$$ 
See Figure~\ref{f_orbit}.

We normalize the associated surgery sequence in the form $f_n:z \to z^2 + v_n$ where each
$v_n$ is real and  $v_n \to v:=-2$ from above, and so that the critical point $0$ is periodic 
for $f_n$ with period $n+2$ consisting of real numbers satisfying
$$-2 < v_n = f_n(0) < 0 = f_n^{n+2}(0) < f_n^{n+1}(0) < \cdots < f_n^2(0) < 2$$
For instance, $f_1$ has critical orbit 
$$0 \to v_1 \to v_1^2 + v_1 \to v_1^4+2v_1^3 + v_1^2 + v_1 = 0$$
so that $v_1 \sim -1.7549$. Likewise we can compute $v_2 \sim -1.9408$, $v_3 \sim -1.9854$ and
so on. The first three values for $(v_n - v)\cdot \mu^n$ are approximately
$0.9805$, $0.9472$, $0.9328$; other approximate values are tabulated below.

\begin{table}[ht]\label{my_table}
\begin{tabular}{c|c}
$v_n$ & $(v_n - v)\cdot \mu^n$ \\
\hline
-1.754877666 & 0.980489335 \\
-1.940799806 & 0.947203095 \\
-1.985424253 & 0.932847805 \\
-1.996376137 & 0.927708745 \\
-1.999095682 & 0.926021297 \\
-1.999774048 & 0.925496551 \\
\end{tabular}
\end{table}

We shall see in \S~\ref{Misiurewicz_section} that this sequence converges, and in fact

$$\lim_{n \to \infty} (v_n - v) \cdot \mu^n = \lim_{n\to\infty} (v_n - (-2)) \cdot 4^n = 3\pi^2/32 \sim 0.9252754$$ 
where $3\pi^2/32 = -3/8 \cdot -\pi^2/4$.
We have already seen the factor of $-\pi^2/4$; the factor of $-3/8$ will be explained in
\S~\ref{Misiurewicz_section}.

\section{Surgery sequences for Misiurewicz points}\label{Misiurewicz_section}

In this section we completely analyze surgery sequences in the quadratic family: polynomials of the
form $z \to z^2 + v$ with $v\in \C$. We use the following standard terminology: the {\em Mandelbrot
set} $\M$ is the set of parameters $v$ for which the critical point $0$ of $z \to z^2 + v$ is 
not in the attracting basin of infinity (some authors call this the `filled Mandelbrot set').
For a quadratic polynomial $z \to z^2 + v$
the critical point $\infty$ is periodic, so in order for it to give rise to a surgery sequence the
other critical point $0$ must be strictly preperiodic. Such a parameter $v$ is called a
{\em Misiurewicz point}. 

Let $f:z \to z^2 + v$ be Misiurewicz, and let $O\subset P(0)-0$ be the periodic orbit
with period $m$ and multiplier $\mu$. 
In the sequel we use the following well-known facts about Misiurewicz points:

\begin{proposition}
Let $v$ be a Misiurewicz point. The following are true:
\begin{enumerate}
\item{$v$ is in the Mandelbrot set $\M$;}
\item{the multiplier $\mu$ satisfies $|\mu|>1$; and}
\item{the Julia set $J(f)$ is a dendrite.}
\end{enumerate}
\end{proposition}
For a proof see Douady--Hubbard \cite{Douady_Hubbard}. 

Let $E_\mu$ denote the elliptic curve $\C^*/\langle \mu \rangle$. 
We warn the reader that throughout this section we set up our notation 
slightly differently than in \S~\ref{definition_subsection}.
Choose $p=f^k(0)\in O$ for some least $k>0$ and 
let $U$ be a connected open neighborhood of $p\in O$ 
on which $g:U \to U$ is an attracting local branch of $f^{-m}$. 
Any tail $T$ for $f$ intersects $U$ in a sequence of points which by abuse of notation we denote
$c_{-j}$ for which $g(c_{-j})=c_{-j-1}$, and thus $c_{-n} \to p$ (this indexing and the
definition of $g$ and $U$ differs from \S~\ref{definition_subsection} when $m\ne 1$ but agrees
when $m=1$).

A holomorphic map may be linearized (i.e.\/ holomorphically conjugated to a linear map) 
near an attracting fixed point; see e.g.\/ \cite{Milnor}. Thus 
$g$ is holomorphically conjugate near $p$ to the map $z \to z/\mu$ near $0$.
In particular, $\lim_{n \to \infty} \mu^n(c_{-n}-p)$ exists. We denote by $x(T)$ the image
of this limit in $E_\mu$.
Another way to say this is to use the holomorphic conjugacy to identify $E_\mu$ with the
quotient $(U-p)/\langle g\rangle$ and then the $g$-orbit $\lbrace c_{-n}\rbrace$ becomes identified
with $x(T) \in E_\mu$.

Let $X \subset E_\mu$ denote the closure of the set of points $x(T)$ that arise in this
way.

\begin{proposition}[$X$ is Julia]\label{X_is_Julia}
With notation as above, the set $X$ is equal to $(J(f) \cap (U-p))/\langle g \rangle$ 
under the identification of $(U-p)/\langle g \rangle$ with $E_\mu$.
\end{proposition}
\begin{proof}
The Julia set $J(f)$ is $f$-invariant and therefore also $g$-invariant (where defined). 
Let $Y \subset E_\mu$ denote the image of $(J(f) \cap (U-p))/\langle g \rangle$ under
the identification of $(U-p)/\langle g \rangle$ with $E_\mu$. Since $J(f)$ is $g$-invariant
and closed, it follows that $Y$ is closed.
Since $O$ is a repelling periodic orbit for $f$ contained in $P(0)$, 
it follows that $0\in J(f)$ and therefore also $T\subset J(f)$ for every $T$. Thus
$x(T)\in Y$ so that $X\subset Y$. Conversely, the complete backward orbit of every element
of $J(f)$ is dense in $J(f)$, so that every point of $J(f)$ may be approximated by
some element of some tail. Thus $Y \subset X$ so $X=Y$. 
\end{proof}

Let us now construct the surgery sequence associated to a tail $T$ and a Misiurewicz point $v$.
Repelling periodic orbits are structurally stable, so that for all $w\in \C$ sufficiently close
to $v$ there is a unique repelling point $p(w)$ of period $m$
for $f_w:z \to z^2 + w$ close to $p$; furthermore, $p(w)$ depends holomorphically on $w$ with
$p(v)=p$. For $|w-v|$ sufficiently small, there is a unique local branch $g_w$ of $f_w^{-m}$
with $g_w:U \to U$ fixing $p(w)$.

In a similar manner we may define $c_{-j}(w)$ to be the preimage of $0$ under a suitable 
power of $f_w$ close to $c_{-j}$, so that each $c_{-j}(w)$ depends holomorphically on $w$, 
and $g_w(c_{-j}(w)) = c_{-1-j}(w)$ and therefore $c_{-n}(w) \to p(w)$.

On the other hand, we may also define $q_w:=f_w^k(0)$, so that $p_v=q_v=p$. Define
$\nu:=d/dw|_{w=v}(q_w-p_w)$. 

\begin{lemma}\label{nu_nonzero}
With notation as above, $\nu \ne 0$. 
\end{lemma}
\begin{proof}
This is the only non-elementary point in the paper; it follows from Thurston's theorem
\cite{Douady_Hubbard_Thurston} on the uniqueness of pcf maps of a given topological type. 
Let us see how.

Thurston's theorem (\cite{Douady_Hubbard_Thurston}) gives necessary and sufficient
conditions (not relevant here) that a critically finite branched
map from $S^2$ to $S^2$ is equivalent to a rational map, and furthermore that such
a rational map is unique (up to holomorphic conjugacy) provided the map has `hyperbolic
orbifold'. This condition is rather technical to state completely, but we remark that it
is satisfied automatically when $|P(f)|>4$ and thus in our context there are only finitely
many exceptional cases where it does not hold where $\nu \ne 0$ may be checked (numerically)
by hand.

Since $q_w$ and $p_w$ both depend holomorphically on $w$, we may write $q_w = p_w + h(w-v)$ for
some holomorphic function $h$ with $h(0)=0$. The first observation is that $h$ is not
identically zero. For, if it were, the $f_w$ would all be topologically equivalent pcf
maps, and therefore (by Thurston's theorem) holomorphically conjugate. But distinct
elements of the quadratic family are never holomorphically conjugate; the first claim follows.

If $h(z)$ is not identically zero then we can write $h(z)= \alpha z^k + O(z^{k+1})$ for some
$\alpha \ne 0$. If $k>1$ then there is a real $\epsilon>0$ so that if $S$ is the circle
$|z|=\epsilon$, the image of $S$ under $h$ has winding number $k$ around $0$. Choose $n$
sufficiently large so that for all $w\in S+v$, the difference 
$|c_{-n}(w)-p(w)|$ is small compared to the minimum of $|h|$ on $S$. Then as $w-v$ winds
around $S$ the difference $q_w-c_{-n}(w)$ also winds $k$ times around $0$, and therefore
there are $k$ distinct values of $w$ near $v$ for which $q_w = c_{-n}(w)$. But
then for each of these $w$ the map $f_w$ is pcf (indeed $0$ is periodic)
and furthermore these maps are all topologically equivalent, violating Thurston's theorem.

It follows that $k=1$ so that $h'(0):=\nu$ is nonzero, as claimed.
\end{proof}

The winding number argument of Lemma~\ref{nu_nonzero} actually shows for all 
sufficiently large $n$ that there is a unique $v_n$ near $v$ with $q_{v_n} = c_{-n}(v_n)$. 
Evidently the $f_{v_n}$ are the surgery sequence associated to the given tail $T$.

\begin{proposition}\label{vn_limit_exists}
With notation as above, $\lim_{n \to \infty} (v_n -v)\mu^n$ exists, and its image
in $E_\mu$ is equal to $\nu^{-1} x(T)$.
\end{proposition}
\begin{proof}
Fix some small positive $\epsilon>0$ and fix some large $j$ so that 
$c_{-j}(w)$ is contained in $U$ for $|w-v|<\epsilon$. If $\epsilon$ is small enough, then
for each $w$ with $|w-v|<\epsilon$ the multiplier $\mu(w)$ of $f_w^m$ at $p(w)$ satisfies
$|\mu(w)|>1$ and there are a family of maps $\phi_w:U \to \C$ (depending holomorphically on $w$)
with $\phi_w(p(w))=0$ and
$\phi_w(c_{-j}(w))=1$ conjugating $g_w:U \to U$ to the map $z \to \mu(w)^{-1}z$ on 
$\phi_w(U)$. We may suppose for concreteness that the disk of radius $2$ about $0$ is
contained in all $\phi_w(U)$ for $|w-v|<\epsilon$.

For $z:=w-v$ with $|z|<\epsilon$ let $f(z)=\phi_w(q(w))$. Then $f(0)=0$ and
$\rho:=f'(0) = \nu \phi_v'(p)$ and we may write $f(z)=z\rho/(1+zb_1+z^2b_2 + \cdots)$
and $h(z):=\mu(w)^{-1} = \mu^{-1}(1 + za_1 + z^2a_2 + \cdots)$ for power series uniformly convergent
on $|z|<\epsilon$. With this notation, $v_n=v+z_n$ where $z_n$ is the solution to 
$f(z_n) =h(z_n)^n$.

An elementary estimate (Lemma~\ref{linear_computation}) shows that 
$\lim_{n \to \infty} z_n\mu^n = \rho^{-1}$. But
$\rho^{-1} = \nu^{-1} (\phi_v'(p))^{-1}$ and if we choose $j$ large enough so that $|c_{-j}-p|$ is
small, then $(\phi_v'(p))^{-1}$ is approximately equal to $(c_{-j}-p)^{-1}$.
The claim follows, modulo the proof of Lemma~\ref{linear_computation}.
\end{proof}

We now prove the desired estimate, completing the proof of Proposition~\ref{vn_limit_exists}.

\begin{lemma}\label{linear_computation}
Let $\mu,\rho \in \C$ with $0<|\mu^{-1}|<1$ and $\rho \ne 0$. Let $h(z):=\mu^{-1}(1+za_1 + z^2a_2 + \cdots)$ and 
$f(z):=z\rho/(1 + zb_1 + z^2b_2 + \cdots)$ be holomorphic in some open neighborhood of $0$.
Then for $n \gg 1$ there is a unique $z_n\in \C$ with $|z_n|\ll 1$ such that
$f(z_n)=h(z_n)^n$, and furthermore
$\lim_{n \to \infty} z_n\mu^n = \rho^{-1}$.
\end{lemma}
\begin{proof}
By definition $z_n$ is the solution to
$$z\rho = \mu^{-n}(1 + z a_1 + z^2a_2 + \cdots)^n (1+zb_1 + z^2 b_2 + \cdots) := \mu^{-n} \bigl(1+\sum_{j>0} \kappa_{j,n} z^j\bigr)$$
for suitable coefficients $\kappa_{j,n}$ depending on $n$.
Since $f(z)$ and $h(z)$ are holomorphic in an open neighborhood of $0$, there are positive 
real constants $\alpha$ and $\beta$ so that $|a_k|<\alpha^k$ and $|b_k|<\beta^k$ for all
$k$, and therefore we may estimate $|\kappa_{j,n}|\le \frac {(n+1+j)!} {n!j!}\sigma^j$ 
where $\sigma = \max(\alpha,\beta)$.

Set $\zeta_1:=\mu^{-n}\rho^{-1}$ and recursively define
$$\zeta_{k+1}: = \mu^{-n}\rho^{-1}(1+\sum_{j>0}\kappa_{j,n}\zeta_k^j)$$
Then it follows by induction for $n$ sufficiently large 
that there are constants $C_1>0$ and $0<C_2<1$
independent of $n$ and $k$ so that 
$$|\zeta_{k+1} - \zeta_k| < C_1 \mu^{n(k+1)C_2}$$
(in fact we can take $C_2$ to be any fixed positive number $<1$ at the cost of adjusting $C_1$).
In particular, for each fixed $n$, the $\zeta_k$ converge at a geometric rate to $z_n$, and
by inspection as $n\to \infty$ we have $z_n\mu^n \to \rho^{-1}$.
\end{proof}

\begin{example}
Let's return to the example we worked out in \S~\ref{2_example}. The fixed points of
the quadratic map $f_w:z \to z^2 + w$ are the roots of $z^2 - z + w$, which are
$(1 \pm \sqrt{1-4w})/2$. Thus for our family where $v=-2$ and $|w-v|$ is small, 
the root $p(w)$ is equal to
$(1+\sqrt{1-4w})/2$ and $p'(-2)= -1/3$. On the other hand, $q(w)=w^2+w$ so that $q'(-2)=-3$.
Thus $\nu:=q'(-2)-p'(-2)=-3+1/3 = -8/3$. 
\end{example}

Every $f_w$ in a surgery family is pcf, and therefore $w$ is in the filled Mandelbrot set.
Let $v$ be a Misiurewicz point, and let $\M-v\subset \C$ denote the result of translating the Mandelbrot
set $\M$ so that $v$ is moved to the origin. For each $n$ let $\mu^n(\M-v)$ be the subset
of $\C$ obtained by multiplying $\M-v$ by $\mu^n$. 

Define $\hat{Z}:=\lim\inf \mu^n(\M-v)$ where the limit is taken in the Hausdorff
topology. In other words, a point $p$ is in $\hat{Z}$ if and only if for every infinite
subsequence of dilations $\mu^{n_i}(\M-v)$ there are points $p_i \in \mu^{n_i}(\M_v)$ with
$p_i \to p$. Evidently $\hat{Z}$ is invariant under multiplication by $\mu$, and therefore
$\hat{Z}-0$ covers a closed subset $Z\subset E_\mu$.

\begin{proposition}
The translate $\nu^{-1} Z$ in $E_\mu$ contains $X$.
\end{proposition}
\begin{proof}
This is a direct consequence of 
Proposition~\ref{vn_limit_exists} and Proposition~\ref{X_is_Julia}
\end{proof}

Conversely, we have the following:

\begin{proposition}
Let $v$ be a Misiurewicz point, and suppose $v_{n_i} \in \M$ is a sequence of
points with $v_{n_i} \to v$ so that $(v_{n_i}-v)\mu^{n_i}$ converge to some limit $z$.
Then the image of $\nu^{-1} z$ in $E_\mu$ is in $X$.
\end{proposition}
\begin{proof}
With notation as above, we let $O$ be the periodic orbit of $f_v$, choose
$f^k(0)=p\in O$, and let $U$ be an open neighborhood of $p$ on which $g:U \to U$
is an attracting local branch of $f^{-m}$. Let $V=g(U)$. Since $v_{n_i} \to v$ we
have $f_{v_{n_i}}^k(0) \to p$ and therefore there is some least $k_i$ for which
$x_i:=f_{v_{n_i}}^{k_i}(0) \in V-U$ and $f_{v_{n_i}}^j(0) \in U$ for $k\le j < k_i$
(in fact, $(k_i - k)/m - n_i$ is a constant).
After passing to a subsequence if necessary, we may assume $x_i \to x\in V-U$.
Since $v_{n_i} \in \M$ it follows that the forward
orbit of $x_i$ under $f_{v_{n_i}}$ is uniformly bounded (independent of $i$) and
therefore the forward orbit of $x$ under $f=f_v$ is uniformly bounded, so that
$x$ is in the filled Julia set of $f$. Since $v$ is Misiurewicz, it follows that
$x$ is in the Julia set of $f$ and therefore that the image of $x$ in $E_\mu$ is
contained in $X$. As in our previous calculations (i.e.\/ Lemma~\ref{linear_computation}), 
this image is equal to $\nu^{-1} z$.
\end{proof}

This concludes our proof of Tan Lei's theorem:

\begin{corollary}[Tan Lei \cite{Lei}, Thm.~5.1]
Let $v\in \M$ be a Misiurewicz point associated to a quadratic polynomial $f:z \to z^2 + v$
with Julia set $J(f)$
and let $p\in P(0)$ be periodic for $f$ with multiplier $\mu$. Let $\nu \ne 0$ be
defined as above. Then the Hausdorff limits
$$\lim_{n \to \infty} \nu^{-1}\mu^n(\M-v) \text{ and } \lim_{n \to \infty} \mu^n(J(f)-p)$$
are equal.
\end{corollary}

\begin{remark}
There is nothing important about the quadratic family so far as 
Proposition~\ref{X_is_Julia} or Proposition~\ref{vn_limit_exists} are concerned, except
that our maps have been constrained to lie in a one (complex) dimensional family.

Let $f$ be an arbitary pcf map with a critical point $c\in C(f)$ and 
repelling periodic orbit $O \subset P(c)-C(f)$.
Let $V$ be a 1-dimensional family of nonconjugate rational maps with $f\in V$. For each critical
point $c\in C(f)$ and each $g\in V$ near $f$ there is a critical point $c(g)\in C(g)$
close to $c$; suppose for all $g\in V$ that $P(c')$ is finite for all $c' \in C(g)-c(g)$.
Choose $p\in O$ equal to $f^k(c)$ for some least $c$. For $g$ near $f$ we can define
$q(g):=g^k(c(g))$ and $p(g)$ to be the repelling periodic point for $g$ near $p$, and
we may define $\nu:= d/dg|_{g=f}q(g)-p(g)$. Then providing $\nu \ne 0$ we may
construct a surgery family $f_n \in V$ associated to any tail $T$ for $O,c$ exactly as
above, and the analog of Proposition~\ref{vn_limit_exists} holds for 
$\lim_{n \to \infty} \mu^n (f_n - f)$. The condition $\nu \ne 0$ holds providing
the pcf maps obtained by deforming $f$ are holomorphically rigid, which according to
Thurston's theorem holds automatically
if $|P(f)|>4$. Presumably the exceptions may be enumerated.

It is worth pointing out that $\nu = 0$ 
if $V$ is a Latt\`es family; see e.g.\/ \cite{Milnor}, and indeed such a family
evidently does not contain a surgery sequence for any $f\in V$.
\end{remark}

\section{Another example}

\begin{figure}[htpb]
\centering
\includegraphics[scale=0.25]{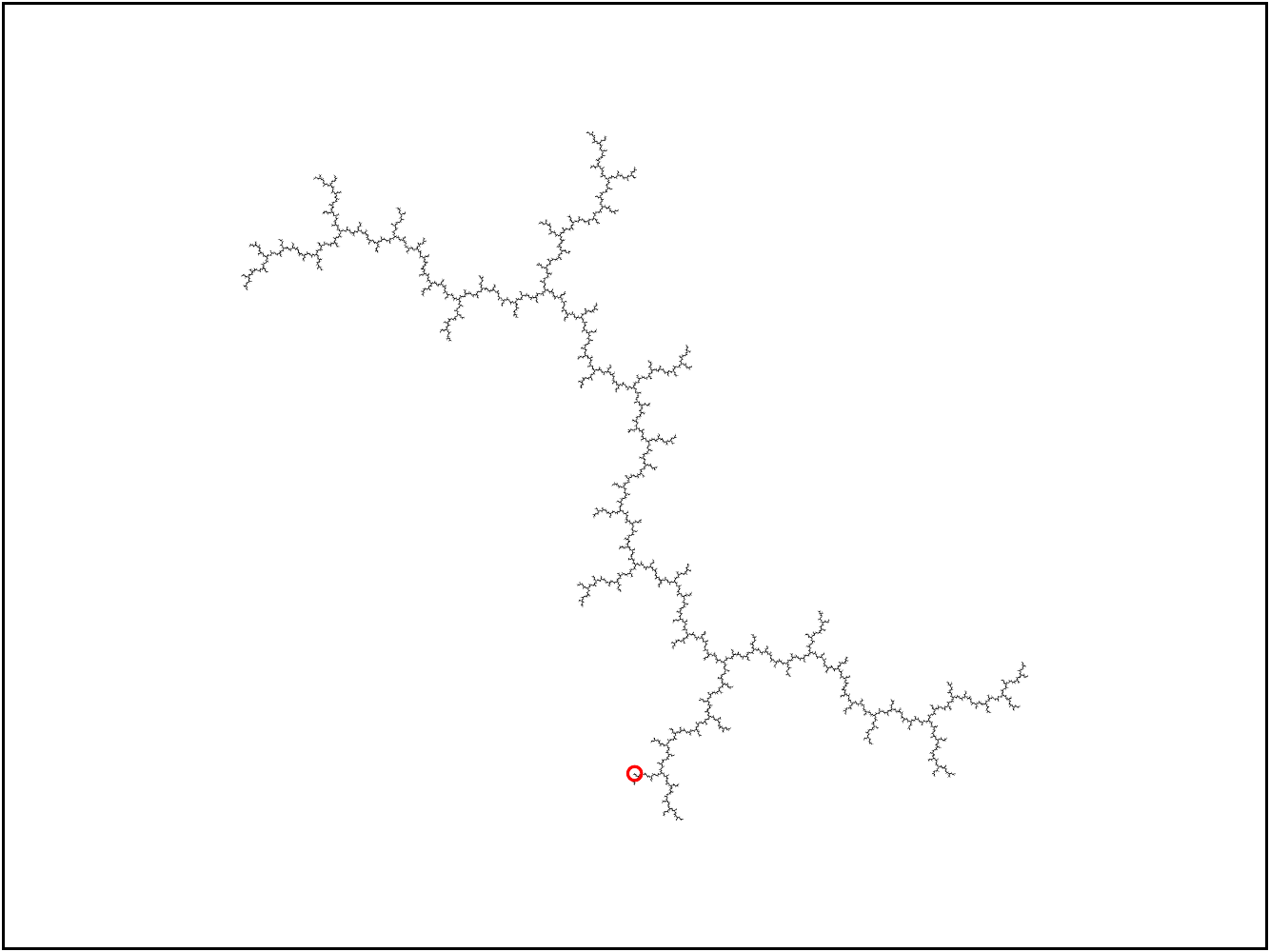} \quad \includegraphics[scale=0.25]{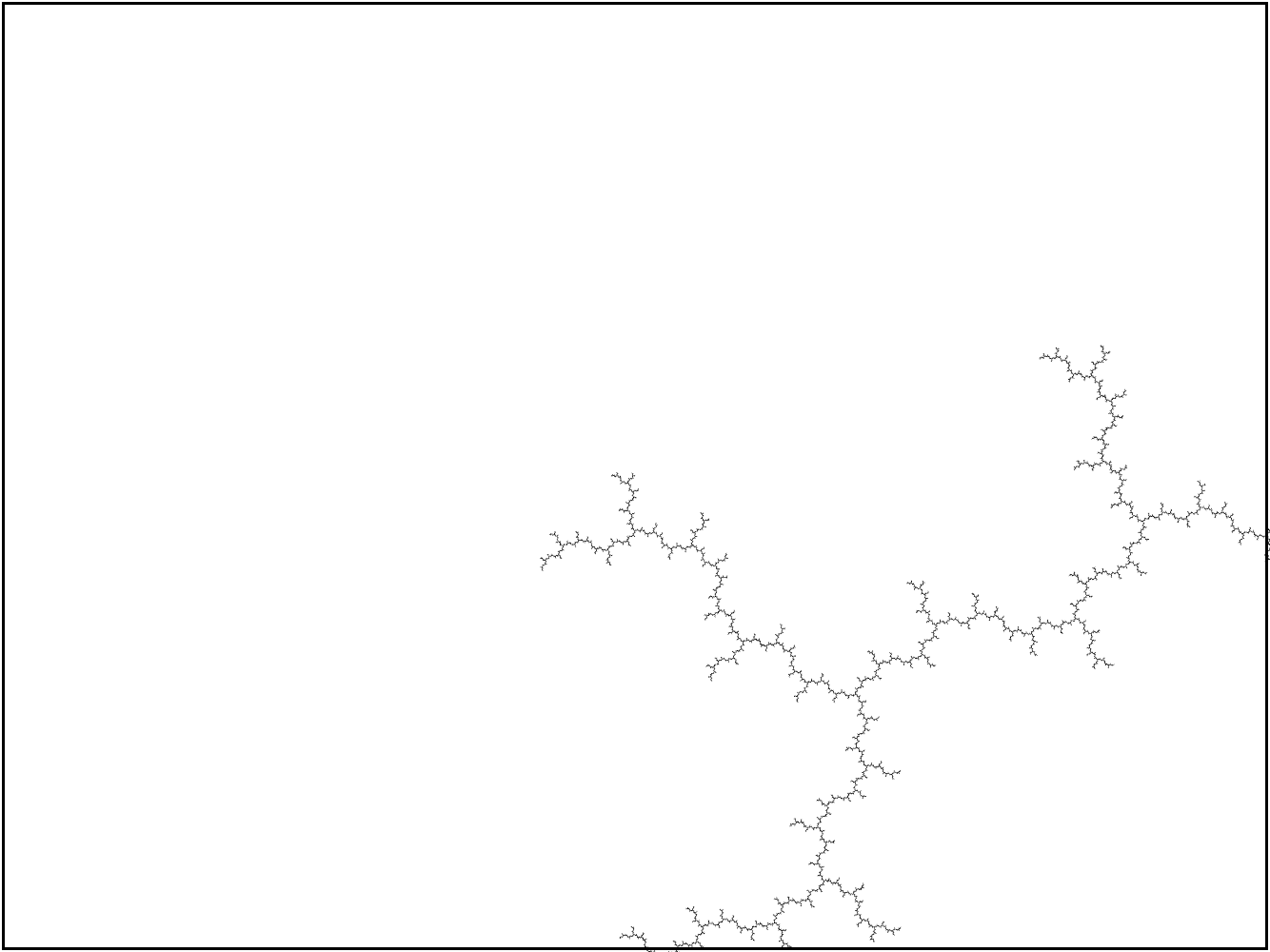}
\caption{The Julia set of $z \to z^2+i$ and a blow-up near $-i$.}
\label{Julia_figure}
\end{figure}

\begin{figure}[htpb]
\centering
\includegraphics[scale=0.25]{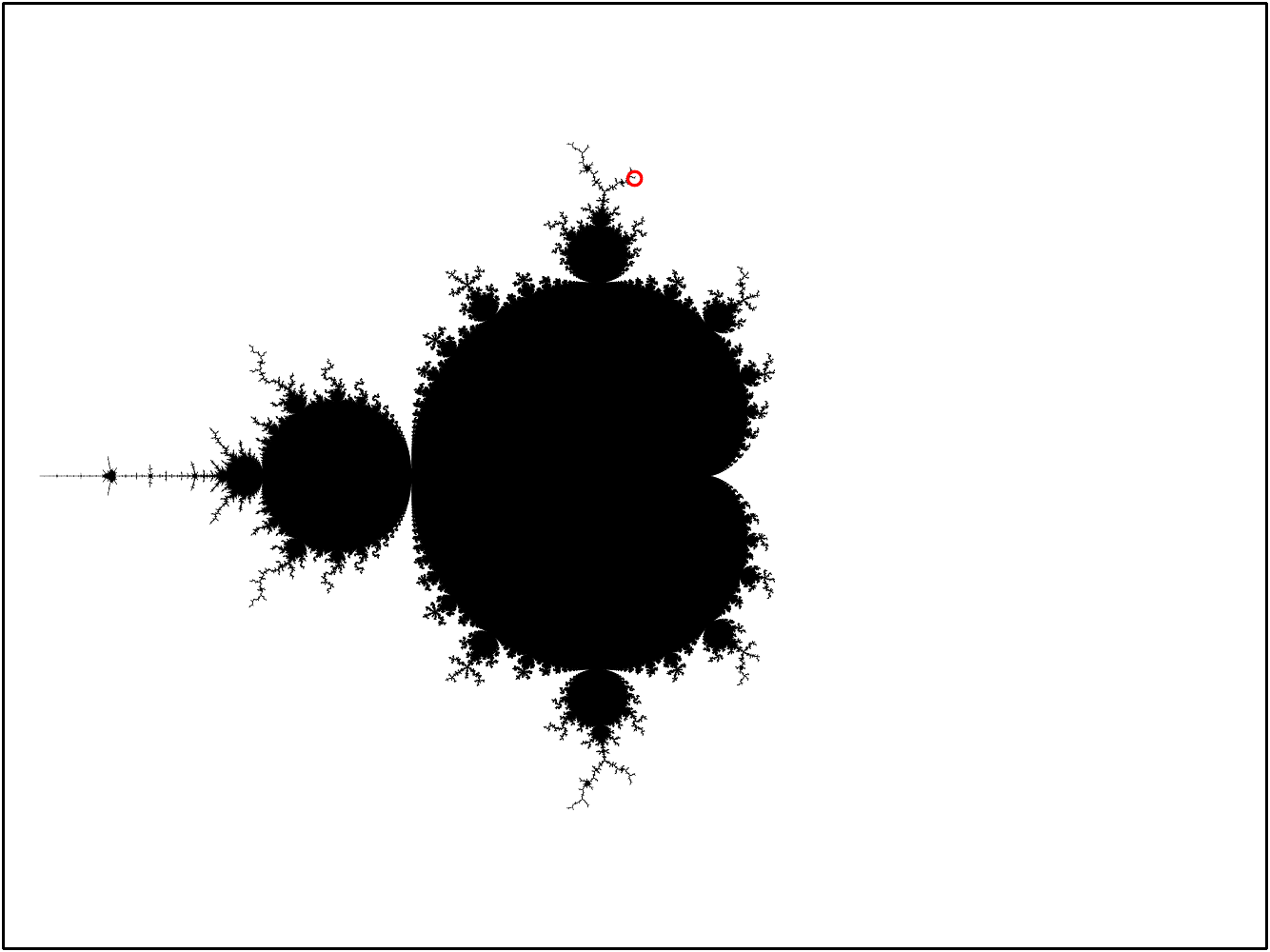} \quad \includegraphics[scale=0.25]{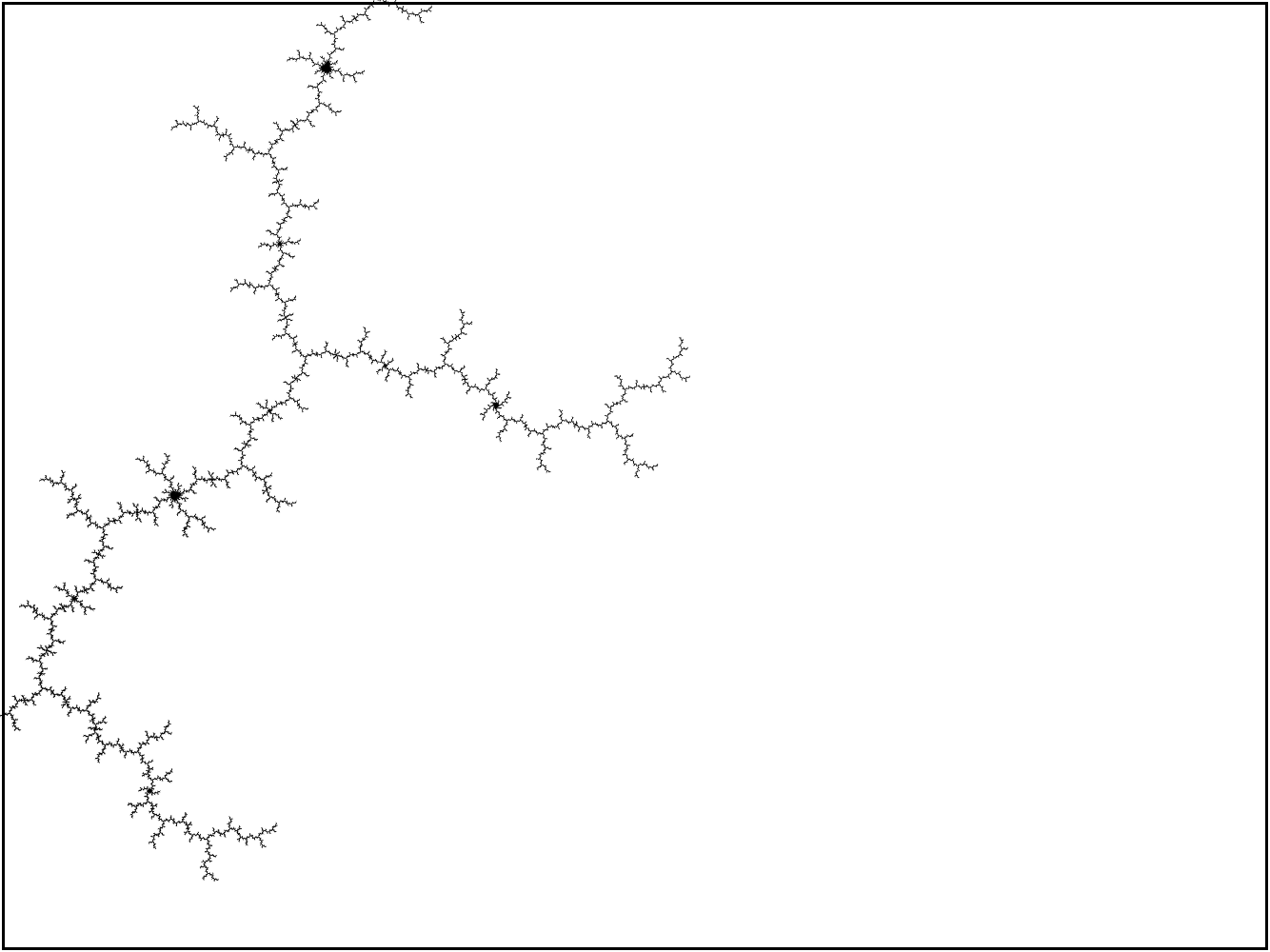}
\caption{The Mandelbrot set and a blow-up near $i$.}
\label{Mandelbrot_figure}
\end{figure}

The map $z: \to z^2 + i$ has $P(0):=\lbrace i, i-1,-i\rbrace$ and we can take
$O:=\lbrace i-1,-i\rbrace$ and $p=-i \in O$. The multiplier is $\mu:=4(1+i)$. The
Julia set $J$ and a blow-up near $-i$ are illustrated in Figure~\ref{Julia_figure}. 
The Mandelbrot set $\M$ and a blow-up near $i$ are illustrated in 
Figure~\ref{Mandelbrot_figure} for comparison.

\section{Acknowledgements}

I would like to thank Frank Calegari, Sarah Koch, Curt McMullen, Alden Walker and the
anonymous referees for valuable comments and encouragement. 

I met Tan Lei in person only once, at the Thurston memorial conference in 2013, although I 
corresponded with her by email from time to time. I have always loved her theorem on
the asymptotic self-similarity of Julia and Mandelbrot sets, and I regret that I cannot share this 
argument with her (which I believe she would have found amusing) but regret her absence much more.

\end{document}